\documentclass[reqno]{amsart}
\usepackage{amssymb}
\usepackage[colorlinks=false,hypertex]{hyperref}
\allowdisplaybreaks[4]
\theoremstyle{plain}
\newtheorem{thm}{Theorem}
\theoremstyle{remark}
\newtheorem{rem}{Remark}
\DeclareMathOperator{\td}{d\mspace{-2mu}}
\date{Completed on Monday 6 October 2008 in Carlisle B, VU Student Village, Australia}
\date{Revised on 17 December 2009 in Tianjin City}
\date{}

\begin{document}

\title[Complete monotonicity of functions involving polygamma functions]
{Uniqueness of nontrivially complete monotonicity for a class of functions involving polygamma functions}

\author[F. Qi]{Feng Qi}
\address[F. Qi]{Department of Mathematics, College of Science, Tianjin Polytechnic University, Tianjin City, 300160, China}
\email{\href{mailto: F. Qi <qifeng618@gmail.com>}{qifeng618@gmail.com}, \href{mailto: F. Qi <qifeng618@hotmail.com>}{qifeng618@hotmail.com}, \href{mailto: F. Qi <qifeng618@qq.com>}{qifeng618@qq.com}}
\urladdr{\url{http://qifeng618.spaces.live.com}}

\author[B.-N. Guo]{Bai-Ni Guo}
\address[B.-N. Guo]{School of Mathematics and Informatics, Henan Polytechnic University, Jiaozuo City, Henan Province, 454010, China}
\email{\href{mailto: B.-N. Guo <bai.ni.guo@gmail.com>}{bai.ni.guo@gmail.com}, \href{mailto: B.-N. Guo <bai.ni.guo@hotmail.com>}{bai.ni.guo@hotmail.com}}
\urladdr{\url{http://guobaini.spaces.live.com}}

\begin{abstract}
For $m,n\in\mathbb{N}$, let $f_{m,n}(x)=\bigr[\psi^{(m)}(x)\bigl]^2+\psi^{(n)}(x)$ on $(0,\infty)$. In the present paper, we prove using two methods that, among all $f_{m,n}(x)$ for $m,n\in\mathbb{N}$, only $f_{1,2}(x)$ is nontrivially completely monotonic on $(0,\infty)$. Accurately, the functions $f_{1,2}(x)$ and $f_{m,2n-1}(x)$ are completely monotonic on $(0,\infty)$, but the functions $f_{m,2n}(x)$ for $(m,n)\ne(1,1)$ are not monotonic and does not keep the same sign on $(0,\infty)$.
\end{abstract}



\thanks{The authors were partially supported by the China Scholarship Council and the Science Foundation of Tianjin Polytechnic University}

\thanks{This paper was typeset using \AmS-\LaTeX}

\maketitle


\section{Introduction}

Recall~\cite[Chapter~XIII]{mpf-1993} and~\cite[Chapter~IV]{widder} that a function $f$ is said to be completely monotonic on an interval $I$ if $f$ has derivatives of all orders on $I$ and
\begin{equation}\label{CM-dfn}
(-1)^{n}f^{(n)}(x)\ge0
\end{equation}
for $x\in I$ and $n\ge0$. The celebrated Bernstein-Widder's Theorem \cite[p.~160, Theorem~12a]{widder} states that a function $f(x)$ is completely monotonic on $[0,\infty)$ if and only if there exists a bounded and non-decreasing function $\alpha(t)$ such that
\begin{equation} \label{berstein-1}
f(x)=\int_0^\infty e^{-xt}\td\alpha(t)
\end{equation}
converges for $x\in[0,\infty)$. This expresses that a completely monotonic function $f$ on $[0,\infty)$ is a Laplace transform of the measure $\alpha$.
\par
It is general knowledge that the classical Euler gamma function
\begin{equation}\label{gamma-dfn}
\Gamma(x)=\int^\infty_0t^{x-1} e^{-t}\td t,\quad x>0,
\end{equation}
the psi function $\psi(x)=\frac{\Gamma'(x)}{\Gamma(x)}$, and the polygamma functions $\psi^{(i)}(x)$ for $i\in\mathbb{N}$ are a series of important special functions, since they have much extensive applications in many branches such as statistics, probability, number theory, theory of $0$-$1$ matrices, graph theory, combinatorics, physics, engineering, and other mathematical sciences.
\par
For $m,n\in\mathbb{N}$, let
\begin{equation}
f_{m,n}(x)=\psi^{(n)}(x)+\bigr[\psi^{(m)}(x)\bigl]^2, \quad x>0.
\end{equation}
The positivity and complete monotonicity of $f_{1,2}(x)$ have been proved, recovered, generalized, and applied in a number of literature in recent years. For much detailed information on these, please refer to \cite{AAM-Qi-09-PolyGamma.tex, notes-best-simple-cpaa.tex} and closely-related references therein.
\par
The aim of this paper is to demonstrate using two methods that, among all $f_{m,n}(x)$ for $m,n\in\mathbb{N}$, only $f_{1,2}(x)$ is nontrivially completely monotonic on $(0,\infty)$.
\par
Our main results may be stated as follows.

\begin{thm}\label{poly-square-poly-thm}
For $m,n\in\mathbb{N}$, the functions $f_{1,2}(x)$ and $f_{m,2n-1}(x)$ are completely monotonic on $(0,\infty)$, but the functions $f_{m,2n}(x)$ for $(m,n)\ne(1,1)$ are not monotonic and does not keep the same sign on $(0,\infty)$.
\end{thm}

\begin{rem}
As shown in the proof of Theorem~\ref{poly-square-poly-thm} below, the complete monotonicity of the functions $f_{m,2n-1}(x)$ is trivial, so $f_{1,2}(x)$ is the only nontrivial completely monotonic function on $(0,\infty)$ among all $f_{m,n}(x)$ for $m,n\in\mathbb{N}$.
\end{rem}

\section{First proof of Theorem~\ref{poly-square-poly-thm}}\label{long-proof-sec}

The following formulas are listed in~\cite{abram}: For $r>0$, $x>0$ and $n\in\mathbb{N}$, we have
\begin{align}
\label{gam}
\frac1{x^r}&=\frac1{\Gamma(r)}\int_0^\infty t^{r-1}e^{-xt}\td t,\\
\label{psin}
\psi ^{(n)}(x)&=(-1)^{n+1}\int_{0}^{\infty}\frac{t^{n}}{1-e^{-t}}e^{-xt}\td t,\\
\label{psisymp4}
\psi^{(n-1)}(x+1)&=\psi^{(n-1)}(x)+\frac{(-1)^{n-1}(n-1)!}{x^n}.
\end{align}
For $n>m\ge1$, direct calculation and utilization of the above three formulas give
\begin{equation*}
f_{m,n}(x)-f_{m,n}(x+1)=\bigr[\psi^{(m)}(x)\bigl]^2-\bigr[\psi^{(m)}(x+1)\bigl]^2 +\psi^{(n)}(x)-\psi^{(n)}(x+1)
\end{equation*}
\begin{align*}
&=\bigr[\psi^{(m)}(x)-\psi^{(m)}(x+1)\bigl]\bigr[\psi^{(m)}(x)+\psi^{(m)}(x+1)\bigl] +\psi^{(n)}(x)-\psi^{(n)}(x+1)\\
&=\frac{(-1)^{m+1}m!}{x^{m+1}}\biggl[2\psi^{(m)}(x)+\frac{(-1)^mm!}{x^{m+1}}\biggr] +\frac{(-1)^{n+1}n!}{x^{n+1}}\\
&=\frac{(-1)^{m+1}2\cdot m!}{x^{m+1}}\biggl[\psi^{(m)}(x)+\frac{(-1)^mm!}{2x^{m+1}} +\frac{(-1)^{n-m}n!}{2\cdot m!x^{n-m}}\biggr]\\
&=\frac{(-1)^{m+1}2\cdot m!}{x^{m+1}}\biggl[(-1)^{m+1}\int_{0}^{\infty}\frac{t^{m}}{1-e^{-t}}e^{-xt}\td t\\ &\quad+\frac{(-1)^m}{2}\int_0^\infty t^{m}e^{-xt}\td t +\frac{(-1)^{n-m}n!}{2\cdot m!(n-m-1)!}\int_0^\infty t^{n-m-1}e^{-xt}\td t\biggr]\\
&=\frac{2\cdot m!}{x^{m+1}}\int_0^\infty\biggl[\frac1{1-e^{-t}} +\binom{n-1}m\frac{(-1)^{n-1}nt^{n-2m-1}}{2}-\frac12\biggr]t^me^{-xt}\td t\\
&=\frac{2\cdot m!}{x^{m+1}}\int_0^\infty\biggl[\frac1{t^{n-2m-1}}\biggl(\frac1{1-e^{-t}}-\frac12\biggr) -(-1)^{n}\frac{n}2\binom{n-1}m\biggr]t^{n-m-1}e^{-xt}\td t\\
&\triangleq \frac{2\cdot m!}{x^{m+1}}\int_0^\infty\biggl[h_{n-2m-1}(t) -(-1)^{n}\frac{n}2\binom{n-1}m\biggr] t^{n-m-1}e^{-xt}\td t
\end{align*}
and
\begin{equation*}
h_k'(t) =\frac{\bigl(1-e^{2t}\bigr)\bigl[k-2te^t/\bigl(1-e^{2t}\bigr)\bigr]}{2t^{k+1}\bigl(e^t-1\bigr)^2} \triangleq \frac{\bigl(1-e^{2t}\bigr)[k-\omega(t)]}{2t^{k+1}\bigl(e^t-1\bigr)^2},
\end{equation*}
where $k$ is any given integer and $\omega(t)$ is increasing on $(0,\infty)$ with the limits
$$
\lim_{t\to0^+}\omega(t)=-1\quad \text{and}\quad \lim_{t\to\infty}\omega(t)=0.
$$
Hence, if $k\ge0$ then the function $h_k'(t)$ is negative on $(0,\infty)$, and if $k\le-1$ then the derivative $h_k'(t)$ is positive on $(0,\infty)$. This means that the function $h_k(t)$ is decreasing for $k\ge0$ and increasing for $k\le-1$ on $(0,\infty)$. Furthermore, from
\begin{equation*}
\lim_{t\to0^+}h_k(t)=
\begin{cases}
0,& k\le-2;\\
1,&k=-1;\\
\infty,& k\ge0
\end{cases}
\end{equation*}
and
\begin{equation*}
\lim_{t\to\infty}h_k(t)=
\begin{cases}
\infty,& k\le-1;\\
\dfrac12,&k=0;\\
0,& k\ge1;
\end{cases}
\end{equation*}
it follows that
\begin{enumerate}
\item
if $k\le-2$ or $k\ge1$, then $0<h_k(t)<\infty$;
\item
if $k=-1$, then $1<h_k(t)<\infty$;
\item
if $k=0$, then $\frac12<h_k(t)<\infty$.
\end{enumerate}
Therefore, by a recourse to Bernstein-Widder's Theorem \cite[p.~160, Theorem~12a]{widder}, as mentioned in~\eqref{berstein-1}, and a fact that a product of finite complete monotonic functions is also completely monotonic, it follows that
\begin{enumerate}
\item
if $n=2i-1$ for $i\in\mathbb{N}$, the function
$$
f_{m,n}(x)-f_{m,n}(x+1)=f_{m,2i-1}(x)-f_{m,2i-1}(x+1)
$$
is completely monotonic on $(0,\infty)$;
\item
if $n=2i$ for $i\in\mathbb{N}$, the quantity
\begin{equation*}
\frac{n}2\binom{n-1}m=i\binom{2i-1}m
\begin{cases}
=1,& i=m=1;\\
=0,& 2i-1<m;\\
\ge2,&\text{$i\ge2$ and $2i-1\ge m$};
\end{cases}
\end{equation*}
which implies that the function
$$
f_{m,n}(x)-f_{m,n}(x+1)=f_{m,2i}(x)-f_{m,2i}(x+1)
$$
is completely monotonic on $(0,\infty)$ for $i=m=1$ or $2i-1<m$;
\item
if $i\ge2$ and $2i-1\ge m$, then the function $f_{m,2i}(x)-f_{m,2i}(x+1)$ is not completely monotonic on $(0,\infty)$.
\end{enumerate}
\par
If the function $f_{m,n}(x)-f_{m,n}(x+1)$ is completely monotonic on $(0,\infty)$, then
\begin{equation*}
(-1)^\ell[f_{m,n}(x)-f_{m,n}(x+1)]^{(\ell)} =(-1)^\ell f_{m,n}^{(\ell)}(x) -(-1)^\ell f_{m,n}^{(\ell)}(x+1)\ge0
\end{equation*}
which can be rewritten as
\begin{equation*}
(-1)^\ell f_{m,n}^{(\ell)}(x)\ge(-1)^\ell f_{m,n}^{(\ell)}(x+1)\ge\dotsm\ge (-1)^\ell f_{m,n}^{(\ell)}(x+i)\to0
\end{equation*}
as $i\to\infty$, where $\ell\ge0$. Thus, the function $f_{m,n}(x)$ is also completely monotonic on $(0,\infty)$. As a result, the functions $f_{m,2i-1}(x)$ for $2i-1>m$, $f_{1,2}(x)$, and $f_{m,2i}(x)$ for $2i-1<m$ and $2i>m$ are completely monotonic on $(0,\infty)$. In a word, the functions $f_{1,2}(x)$ and $f_{m,2i-1}(x)$ for $2i-1>m$ are completely monotonic on $(0,\infty)$.
\par
The formula~\eqref{psin} shows that the functions $(-1)^{m+1}\psi^{(m)}(x)$ and $\psi^{(2n-1)}(x)$ for $m,n\in\mathbb{N}$ are completely monotonic on $(0,\infty)$, so the function
$$
\bigr[(-1)^{m+1}\psi^{(m)}(x)\bigl]^2+\psi^{(2n-1)}(x) =\bigr[\psi^{(m)}(x)\bigl]^2+\psi^{(2n-1)}(x)=f_{m,2n-1}(x)
$$
for $m,n\in\mathbb{N}$, the result of a square and an additive operation of completely monotonic functions, is completely monotonic on $(0,\infty)$.
\par
In~\cite[Theorem~2.1]{Ismail-Muldoon-119}, \cite[Theorem~2.1]{Muldoon-78} and~\cite[Lemma~1.3]{sandor-gamma-2-ITSF.tex}, the function $\psi(x)-\ln x+\frac{\alpha}x$ was proved to be completely monotonic on $(0,\infty)$, i.e.,
\begin{equation}\label{com-psi-ineq-dfn}
(-1)^i\Bigl[\psi(x)-\ln x+\frac{\alpha}x\Bigr]^{(i)}\ge0
\end{equation}
for $i\ge0$, if and only if $\alpha\ge1$, so is its negative, i.e., the inequality~\eqref{com-psi-ineq-dfn} is reversed, if and only if $\alpha\le\frac12$. In \cite{chen-qi-log-jmaa}, the function $\frac{e^x\Gamma(x)} {x^{x-\alpha}}$ was proved to be logarithmically completely monotonic on $(0,\infty)$, i.e.,
\begin{equation}\label{com-psi-ineq-ln-dfn}
(-1)^k\biggl[\ln\frac{e^x\Gamma(x)} {x^{x-\alpha}}\biggr]^{(k)}\ge0
\end{equation}
for $k\in\mathbb{N}$, if and only if $\alpha\ge1$, so is its reciprocal, i.e., the inequality~\eqref{com-psi-ineq-ln-dfn} is reversed, if and only if $\alpha\le\frac12$. Considering the fact (see, for example, \cite[p.~82]{e-gam-rat-comp-mon}) that a completely monotonic function which is non-identically zero cannot vanish at any point on $(0,\infty)$ and rearranging either~\eqref{com-psi-ineq-dfn} or~\eqref{com-psi-ineq-ln-dfn} leads to the double inequalities: For $x\in(0,\infty)$ and $k\in\mathbb{N}$,
\begin{equation}\label{qi-psi-ineq-1}
\ln x-\frac1x<\psi(x)<\ln x-\frac1{2x}
\end{equation}
and
\begin{equation}\label{qi-psi-ineq}
\frac{(k-1)!}{x^k}+\frac{k!}{2x^{k+1}}< (-1)^{k+1}\psi^{(k)}(x) <\frac{(k-1)!}{x^k}+\frac{k!}{x^{k+1}}.
\end{equation}
Using the inequality~\eqref{qi-psi-ineq} in the equation
\begin{equation*}
f_{m,2n}'(x)=2\psi^{(m)}(x)\psi^{(m+1)}(x)+\psi^{(2n+1)}(x)
\end{equation*}
yields that
\begin{equation*}
f_{m,2n}'(x)\le \frac{p_{m,n}(x)}{4x^{2m+2n+3}}
\end{equation*}
and
\begin{align*}
f_{m,2n}'(x)&\ge \frac{q_{m,n}(x)}{2x^{2m+2n+3}},
\end{align*}
where
\begin{align*}
p_{m,n}(x)&=4(2n)!x^{2m+2}+4(2n+1)!x^{2 m+1} -m!(m+1)!x^{2n}\\
&\quad-2[(m!)^2+(m-1)!(m+1)!]x^{2n+1}-4(m-1)!m!x^{2n+2}
\end{align*}
and
\begin{align*}
q_{m,n}(x)&=2(2n)!x^{2m+2}+(2n+1)!x^{2 m+1}-2m!(m+1)!x^{2n}\\
&\quad-2[(m!)^2+(m-1)!(m+1)!]x^{2n+1}-2(m-1)!m!x^{2n+2}.
\end{align*}
If $m>n$, then the powers of $x$ in $q_{m,n}(x)$ have the relations
\begin{equation}\label{relation-1}
2m+2>2n+2>2n+1>2n\quad \text{and} \quad 2m+2>2m+1,
\end{equation}
so the function $q_{m,n}(x)$ tends to $\infty$ as $x\to\infty$, and so $f_{m,2n}'(x)$ is positive when $x$ is large enough, that is, the function $f_{m,2n}(x)$ is increasing on some infinite interval whose left end point is large enough. If $m<n$, then the powers of $x$ in $q_{m,n}(x)$ have the relations
\begin{equation}\label{relation-2}
2m+1<2n<2n+1<2n+2\quad \text{and}\quad 2m+1<2m+2,
\end{equation}
thus the function $q_{m,n}(x)$ tends to $0^+$ as $x\to0^+$, and so $f_{m,2n}'(x)$ is positive when $x$ is enough close to $0$ from the right-hand side, that is, the function $f_{m,2n}(x)$ is increasing on some interval where the number $0$ is its left end point. If $m=n$, then
\begin{align*}
q_{m,m}(x)&=2[(2m)!-(m-1)!m!]x^{2m+2}+[(2m+1)!-2(m!)^2\\
&\quad-2(m-1)!(m+1)!]x^{2m+1}-2m!(m+1)!x^{2m};
\end{align*}
From $(2m)!-(m-1)!m!\ge1$, it follows as above that the function $f_{m,2m}(x)$ is increasing on some infinite interval whose left end point is large enough. In a word, the function $f_{m,2n}(x)$ is not decreasing on $(0,\infty)$.
\par
Now consider the function $p_{m,n}(x)$. If $m>n$, then the relations in~\eqref{relation-1} still holds and the quantity $2n$ is the smallest power of $x$ in the polynomial $p_{m,n}(x)$; and so, by the similar argument as above, the function $p_{m,n}(x)$ is negative and the function $f_{m,2n}(x)$ is decreasing on some interval whose left end point is $0$. If $m<n$, the relations in~\eqref{relation-2} is also valid and the number $2n+2$ is the greatest power of $x$ in the polynomial $p_{m,n}(x)$, consequently the function $p_{m,n}(x)$ is negative and the function $f_{m,2n}(x)$ is decreasing on some infinite interval whose left end point is large enough. If $m=n$, then
\begin{align*}
p_{m,m}(x)&=4[(2m)!-(m-1)!m!]x^{2m+2}+2[2(2m+1)!-(m!)^2\\
&\quad-(m-1)!(m+1)!]x^{2m+1}-m!(m+1)!x^{2m};
\end{align*}
Therefore, the function $p_{m,m}(x)$ is negative and the function $f_{m,2m}(x)$ is decreasing on some interval whose left end point is $0$. In a word, the function $f_{m,2n}(x)$ is not increasing on $(0,\infty)$.
\par
In conclusion, the functions $f_{m,2n}(x)$ for $(m,n)\ne(1,1)$ are not monotonic on $(0,\infty)$.
\par
By the same argument as discussing the sign of $f_{m,2n}'(x)$ above, we can also turn out that the function $f_{m,2n}(x)$ for $(m,n)\ne(1,1)$ does not keep the same sign on $(0,\infty)$. For the sake of shortening the length of this paper, we omit the detailed argument. The proof of Theorem~\ref{poly-square-poly-thm} is complete.

\section{Second proof of Theorem~\ref{poly-square-poly-thm}}

After going through the proof of Theorem~\ref{poly-square-poly-thm} in Section~\ref{long-proof-sec}, we may outline a sketch of a short proof for Theorem~\ref{poly-square-poly-thm} as follows.
\par
The complete monotonicity of $f_{1,2}(x)$ has been proved in different manners and in more general forms in \cite{AAM-Qi-09-PolyGamma.tex, notes-best-simple-cpaa.tex} and related references therein. For the convenience of the reader, it would have been a good idea to include
a short proof for it, to avoid that the reader has to study either or both of the whole papers \cite{AAM-Qi-09-PolyGamma.tex, notes-best-simple-cpaa.tex}.
\par
Utilizing the formulas~\eqref{gam}, \eqref{psin} and \eqref{psisymp4} to compute straightforwardly gives
\begin{align*}
f_{1,2}(x)-f_{1,2}(x+1)&=\bigl[\psi'(x)-\psi'(x+1)\bigr]\bigl[\psi'(x)+\psi'(x+1)\bigr] \\
&\quad+\bigl[\psi''(x)-\psi''(x+1)\bigr]\\
&=\frac1{x^2}\biggl[2\psi'(x)-\frac1{x^2}\biggr]-\frac{2}{x^3}\\
&=\frac2{x^2}\biggl[\psi'(x)-\frac1{2x^2}-\frac1{x}\biggr]\\
&=\frac2{x^2}\int_0^\infty\biggl(\frac{t}{1-e^{-t}}-1-\frac{t}2\biggr)e^{-xt}\td t\\
&=\frac2{x^2}\int_0^\infty\biggl(\frac{t/2}{\tanh(t/2)}-1\biggr)e^{-xt}\td t,
\end{align*}
which shows that the function discussed is completely monotonic on $(0,\infty)$. So is then the sum
$$
\sum_{k=0}^N[f_{1,2}(x+k)-f_{1,2}(x+k+1)]=f_{1,2}(x)-f_{1,2}(x+N+1)
$$
which converges to $f_{1,2}(x)$ pointwise as $N\to\infty$, and therefore $f_{1,2}(x)$ is completely monotonic on $(0,\infty)$.
\par
Since $(-1)^{n+1}\psi^{(n)}(x)$ is a completely monotonic function for all $n\in\mathbb{N}$, adding the completely monotonic function $\psi^{(2n-1)}(x)$ to the completely monotonic function $\bigl[\psi^{(m)}(x)\bigr]^2$ leads to the complete monotonicity of $f_{m,2n-1}(x)$ for $m,n\in\mathbb{N}$ on $(0,\infty)$.
\par
Using the substitution $xt=u$ in \eqref{psin} yields
\begin{equation}\label{psi-n-der-infinity}
  \psi^{(n)}(x)\sim(-1)^{n+1}\frac{(n-1)!}{x^n}, \quad x\to\infty
\end{equation}
and
\begin{equation}\label{psi-n-der-0}
  \psi^{(n)}(x)\sim(-1)^{n+1}\frac{n!}{x^{n+1}},\quad x\to0^+.
\end{equation}
From these, it follows that
\begin{enumerate}
\item
for large $x$ the behaviour of $f_{m,2n}(x)$ is like
$$
\frac{[(m-1)!]^2}{x^{2m}}\biggl\{1-\frac{(2n-1)!}{[(m-1)!]^2}x^{2(m-n)}\biggr\}
$$
and this is clearly
\begin{enumerate}
\item
negative when $m>n$ and $x$ is large, or
\item
positive when $m<n$ and $x$ is large,
\end{enumerate}
and if $m=n$ the sign is determined by
$$
1-m\binom{2m-1}{m-1}
$$
which is negative unless $m=1$,
\item
for small $x$ the behaviour of $f_{m,2n}(x)$ is like
\begin{equation*}
  \frac{(m!)^2}{x^{2(m+1)}}\biggl\{1-\frac{(2n)!}{(m!)^2}x^{2(m-n)+1}\biggr\}
\end{equation*}
and this is clearly
\begin{enumerate}
\item
negative when $m<n$ and $x$ is small, or
\item
positive when $m\ge n$ and $x$ is small.
\end{enumerate}
\end{enumerate}
As a result, the function $f_{m,2n}(x)$ for $(m,n)\ne(1,1)$ or its negative is not completely monotonic and does not keep the same sign on $(0,\infty)$. The proof of Theorem~\ref{poly-square-poly-thm} is complete.

\subsection*{Acknowledgements}
The original version of this manuscript was completed during the first author's visit to the RGMIA, Victoria University, Australia between March 2008 and February 2009. The first author would like to express many thanks to Professors Pietro Cerone and Server S.~Dragomir and other local colleagues for their invitation and hospitality throughout this period.

\end{document}